\documentclass[10pt]{article}

\usepackage{graphicx}			
\usepackage{amssymb}
\usepackage{amsmath}

\usepackage{verbatim}

\parskip = \baselineskip

\def\vp{\vspace{\baselineskip}}

\title{The Pegasus Tiles: an aperiodic pair}
\date{}
\author{Chaim Goodman-Strauss\\ Univ. Arkansas\\ \tt strauss@uark.edu}

\begin{document}

\maketitle

Here we formally present an aperiodic pair of tiles, with tip-to-tip matching rules, we will call the ``Pegasus" tiles, so named to reflect their curious history. These tiles are mutually locally decomposible to the ``$1+\epsilon+\epsilon^2$'' tiles introduced by Penrose in~\cite{penroseEpsilon}, and are closely related to the Taylor-Socolar aperiodic monotile of~\cite{taylorSocolar}. 

\centerline{\includegraphics[angle=90]{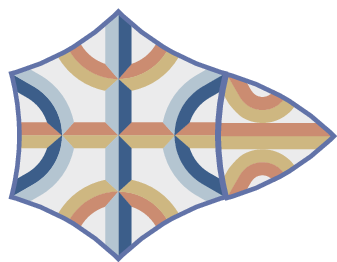}}

In early 1995, before that work appeared, shortly before a meeting on quasicrystals at Les Houches, for an hour or two I believed the tile shown below was aperiodic. As we can see, it is not.\footnote{The tile has ``isohedral number" 3  --- in later years, much larger isohedral number examples have been found: see~\cite{myers} for examples and~\cite{gs_tassellazioni} for further discussion.} Note that the tile is almost exactly the same as the one found by Taylor, but missing some needed structure.

\centerline{\includegraphics[scale=.4]{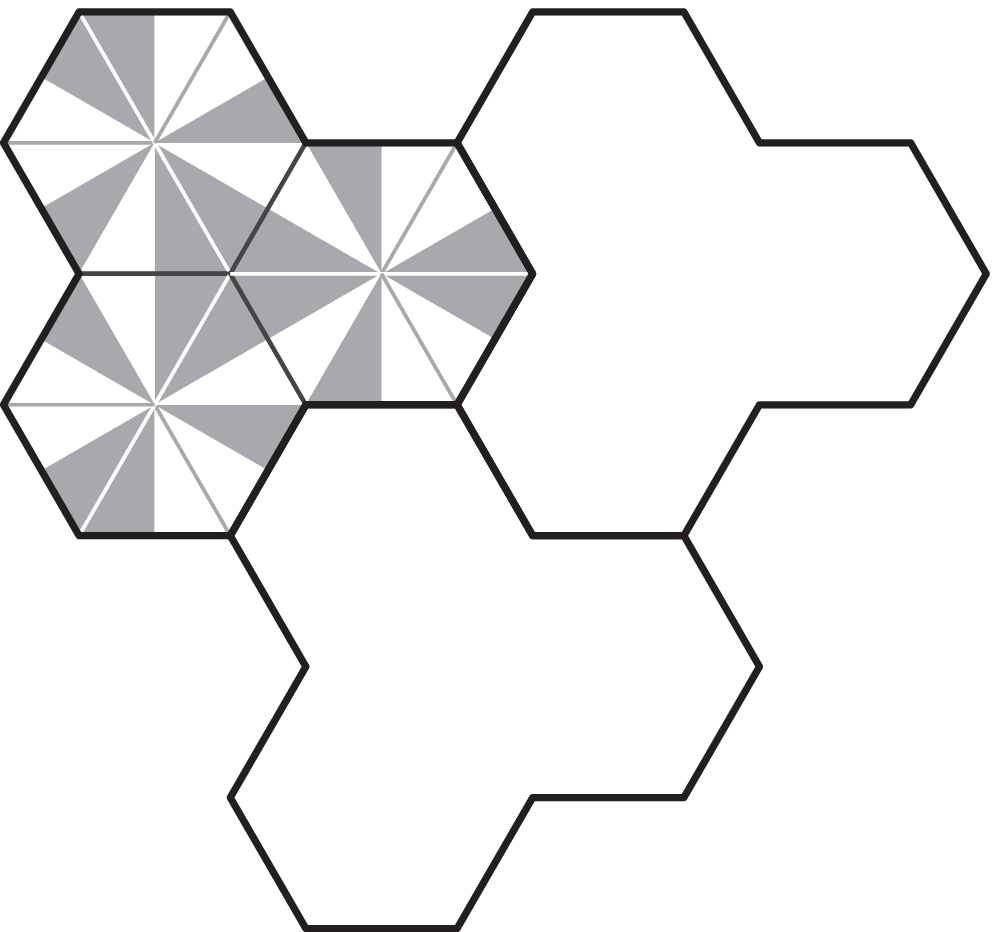}}

I arrived a little late to the meeting, and there was a bit of a buzz in the air: Socolar had found an aperiodic monotile! In fact, as it happened, Socolar had found, essentially, the exact same tile, at almost the exact same time, but had not yet noticed the error. Over the next few days, he modified the construction, ultimately producing a set of about a dozen hexagonal tiles that enforced the structure we had hoped for.  

In the summer of 1995, at a meeting at the Fields Institute, Penrose introduced an ``almost" monotile --- a set of three tiles, but the area of one could be made arbitrarily smaller than the area of the second, which could be made arbitrarily smaller than the area of the third~\cite{penroseEpsilon}. In fact, this phenomenon is relatively common: the tiniest tiles are ``vertex tiles", the middle ones ``edge tiles" and the large ones correspond to the interiors of some original set of tiles: the constructions in, say,~\cite{gs_small} and~\cite{gs_mrst} are also ``$1+\epsilon+\epsilon^2$'', as are many others.

In summer of 1996, I saw this construction for the first time, and, to my astonishment, this was exactly the structure that Socolar and I had thought we were enforcing earlier! Penrose's triple was mutually locally decomposable to Socolar's dozen or so. 

At that time, I had already used ``tip-to-tip" matching rules, in the Trilobite and Cross aperiodic pair of tiles~\cite{gs_small} and realized that Penrose's construction could be recomposed into an aperiodic pair of tiles, but with such rules. Hence the Pegasus pair, shown here  (also, I hope, making the structure more transparent.) The name reflects Penrose, Socolar and myself, in an colorful way.

In any case:

\vp\vp{\bf Theorem:} {\em The Pegasus pair of tiles is {\em aperiodic} --- that is, they do admit tilings, but only admit non-periodic tilings.}

\vp {\bf Proof:} The proof is standard: we show that the tiles can only cluster into larger combinatorial versions of themselves, and so can only form hierarchical tilings, with a unique hierarchy. Hence the tiles {\em can} tile the plane, since they can tile arbitrarily large patches, and the tiles can {\em only} tile non-periodically, since any translation is smaller than some large pegasus tile but any tiling has unique hierarchy.  

The proof that the tiles can only form larger combinatorial versions of themselves requires brute force case checking which we examine in a depth first manner. To streamline this, we adopt some self-evident graphical conventions. A black dot indicates the next location determining further subcases; a grey dot shows the location of black dots for higher level cases. Red markings indicate contradictions. We begin with:

{\includegraphics{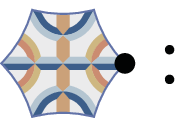}}

\vp
{giving three subcases:} 

\hspace{2em} {\includegraphics{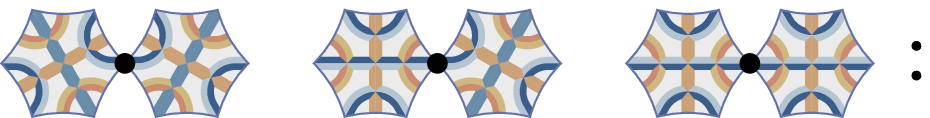}}

\vp The first subcase:
\hspace{2em}\hspace{2em}{\includegraphics{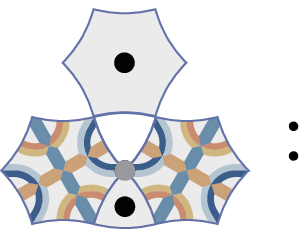}}

\vp (lower)

\centerline{\includegraphics{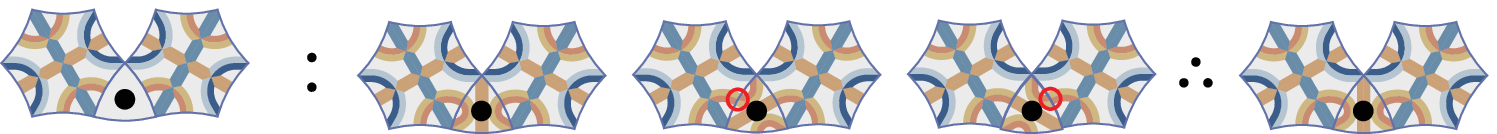}}

\vp (upper)

\centerline{\includegraphics{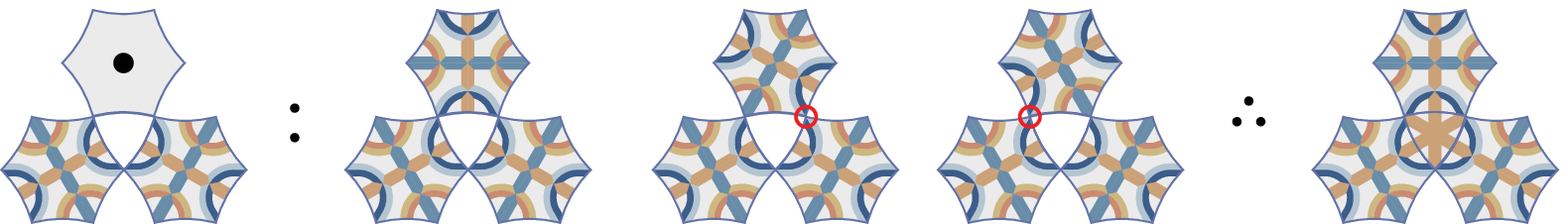}}

\vp Therefore, in the first case

\centerline{\includegraphics{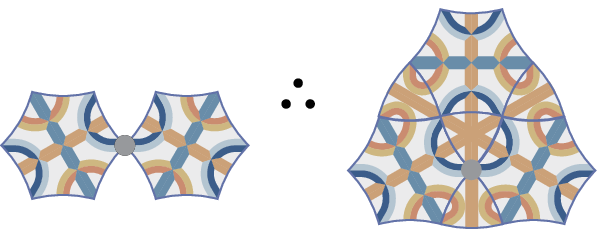}}

\vp\vp 
The second case is more complex:

\vp
{\includegraphics{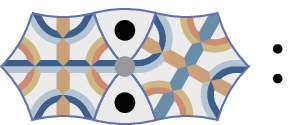}}

\vp We have three possibilities for the upper tile, two of which lead to contradictions when the lower tile is considered. 

\centerline{\includegraphics{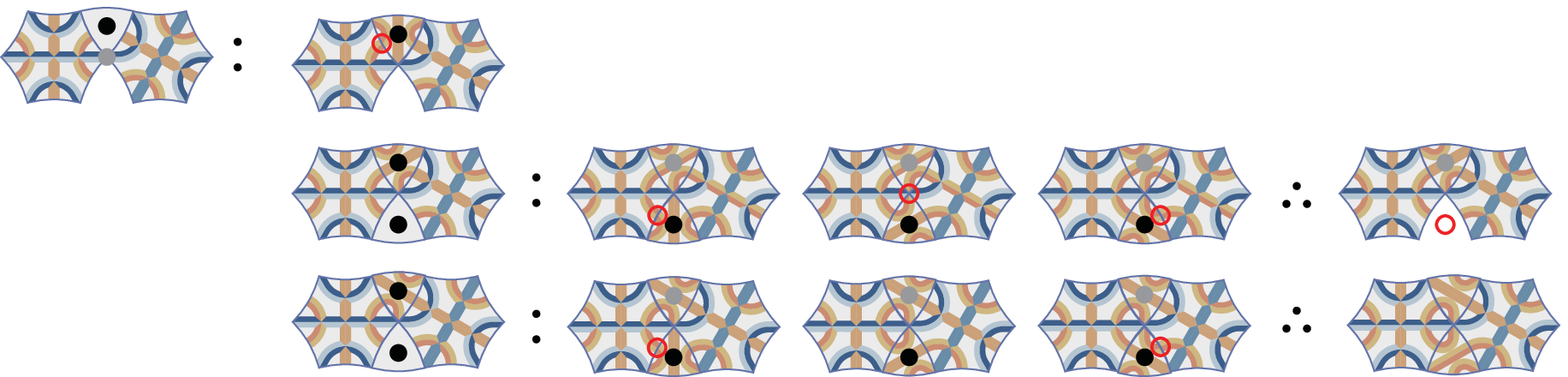}}

\vp Therefore, in the second case, we have 

 {\includegraphics{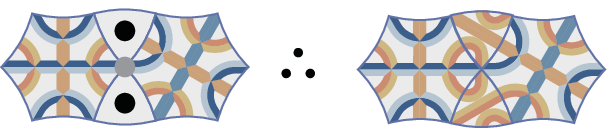}}

\vp 
Next, in the second case, examining nearby hexagonal tiles:

{\includegraphics{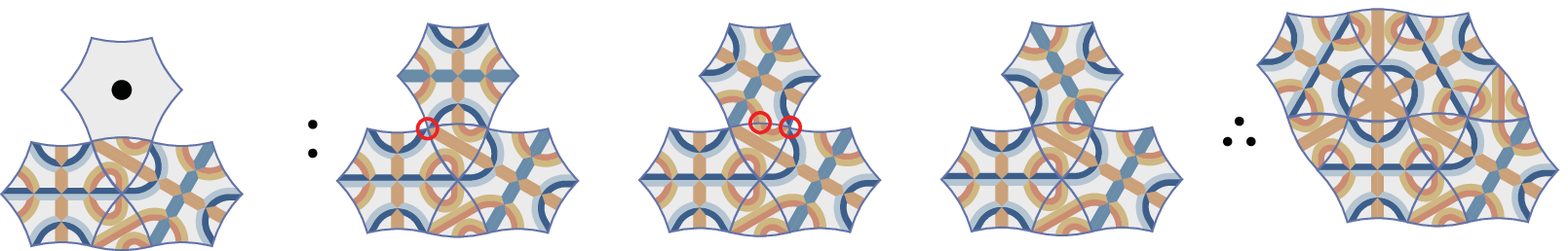}}

\vp 
Thus, the second case can only occur as:

{\includegraphics{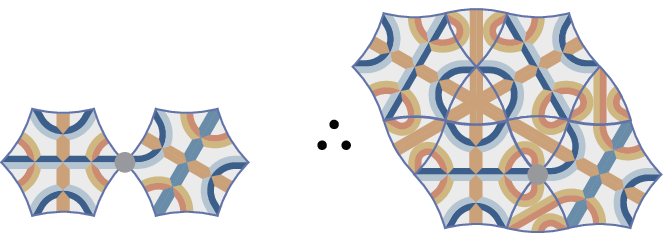}}

\vp\vp We turn to the third case: \hspace{2em} {\includegraphics{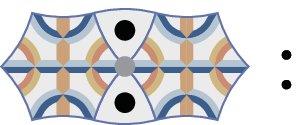}}

\vp Taking into account left-right symmetry, the matching rules allow only one possibility for the upper tile; of the three further possibilities for the lower tile, only one of these satisfies our matching rules.

{\includegraphics{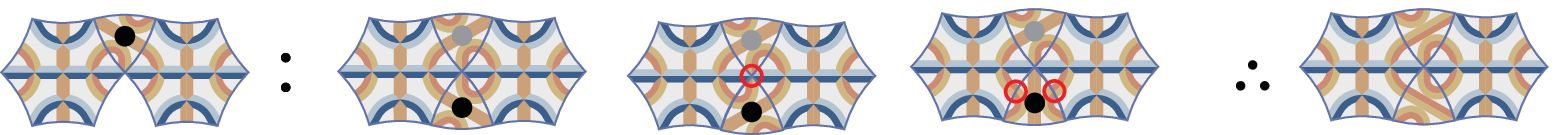}}

\vp Thus:

{\includegraphics{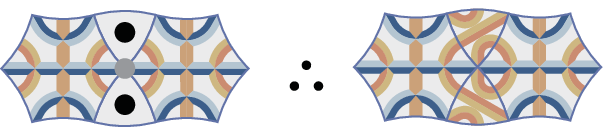}}

\vp Examining possibilities for a neighboring hexagon, 

{\includegraphics{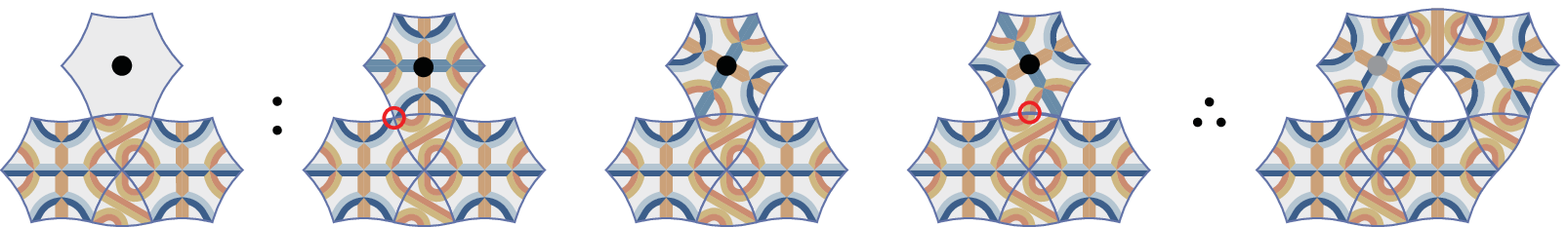}}

\vp 

In short:

{\includegraphics{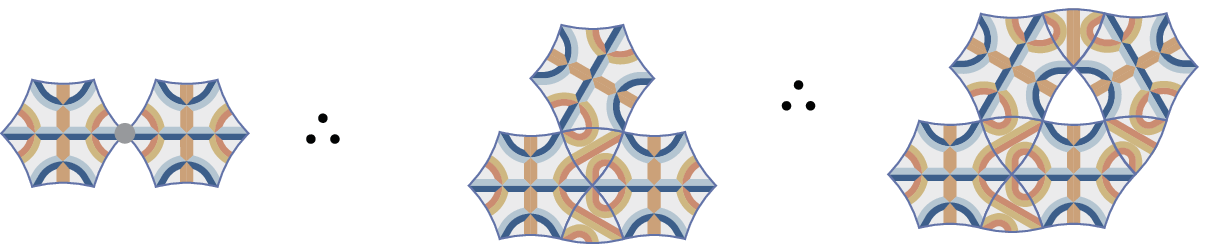}}

\vp 

Of the three possibilities for placing the missing triangle, only one satisfies our matching rules:

{\includegraphics{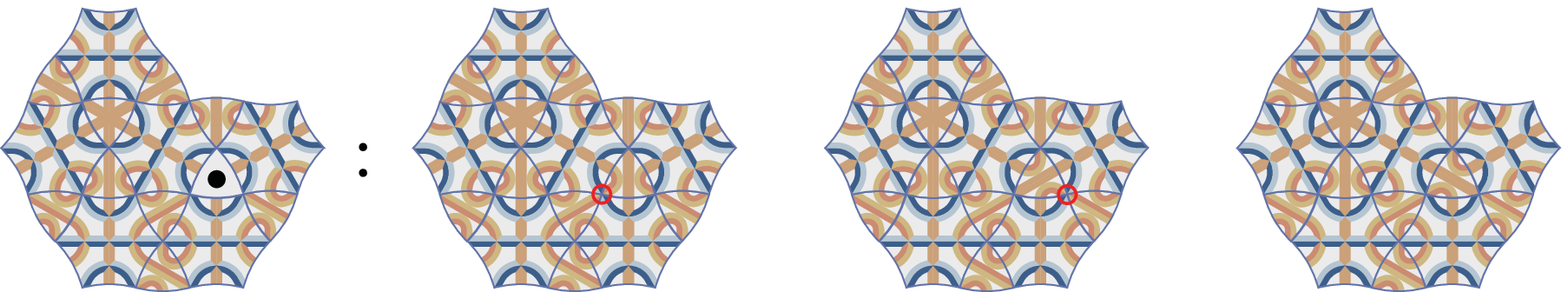}}

\vp 

Our third case can only appear as:

{\includegraphics{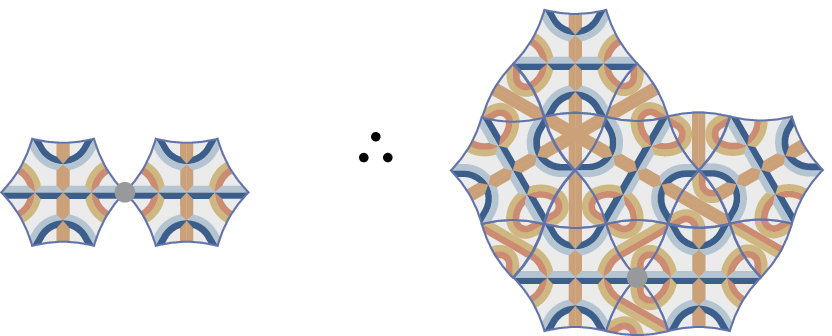}}

\vp\vp

Now we examine the neighborhood of an arbitrary hexagonal tile, in cases two and three. 
 If the gray tile is in case three, it is also in case two with another , and we have

\centerline{\includegraphics{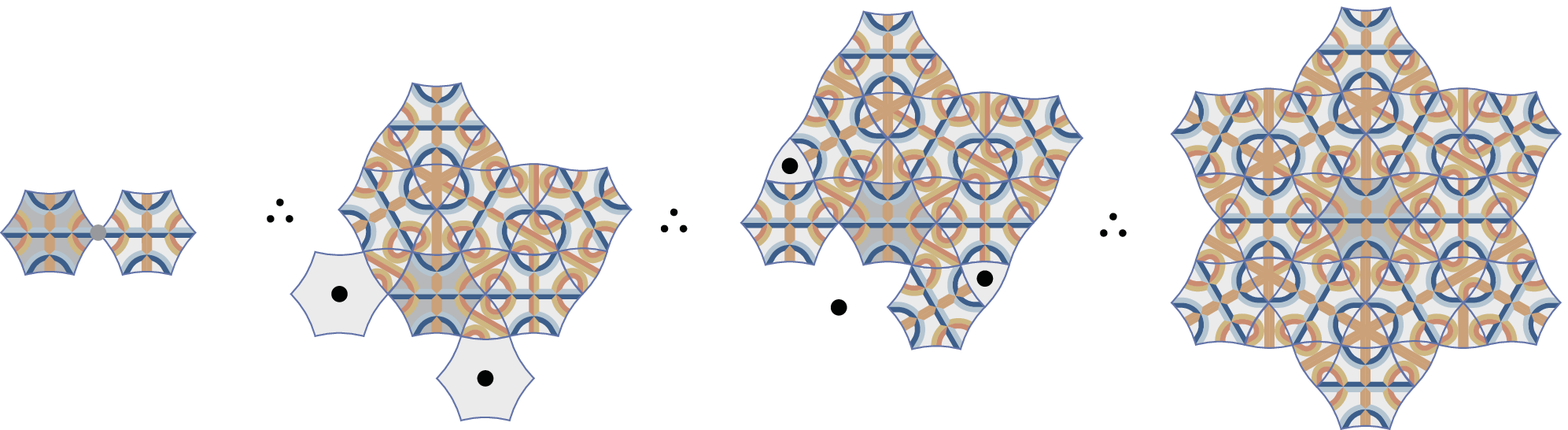}}

In case two, the gray tile is also in case three with another, and so we have

{\includegraphics{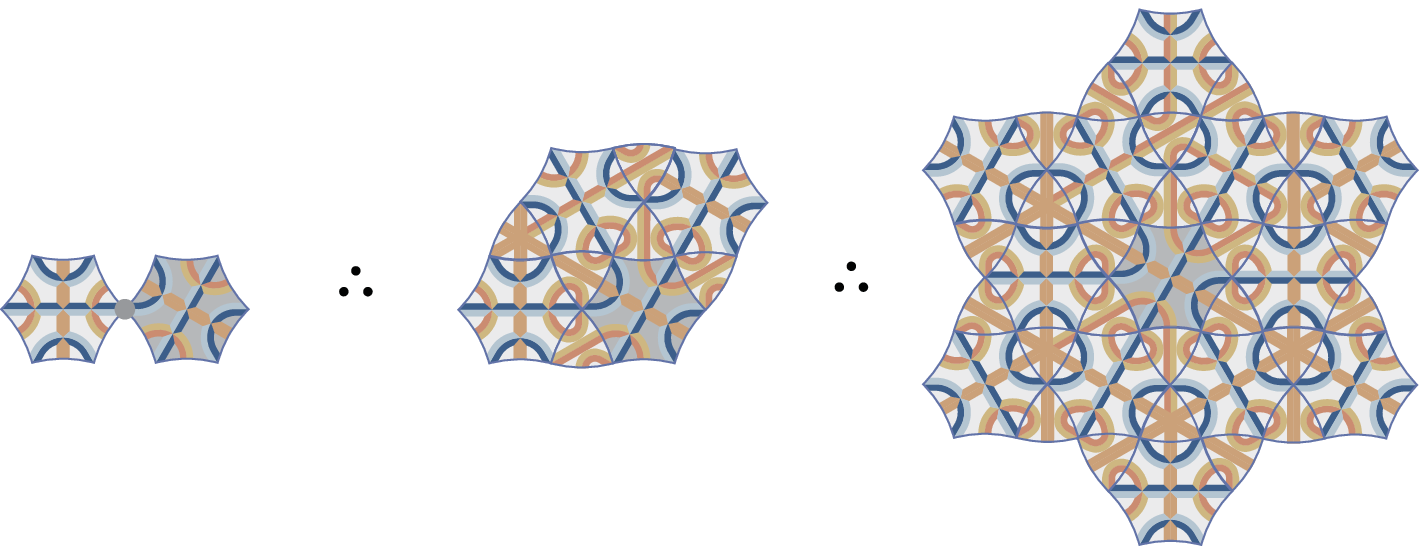}}

\vp
Notice that every tile in case one is also in case 2 with some other tile. 

Therefore, every hexagonal tile is either the center of a large hexagon (cases two and three) or a large triangle (case one) and we may recompose into larger hexagons and triangles, that themselves are combinatorially equivalent to our original tiles, and the proof is complete. 

{\includegraphics{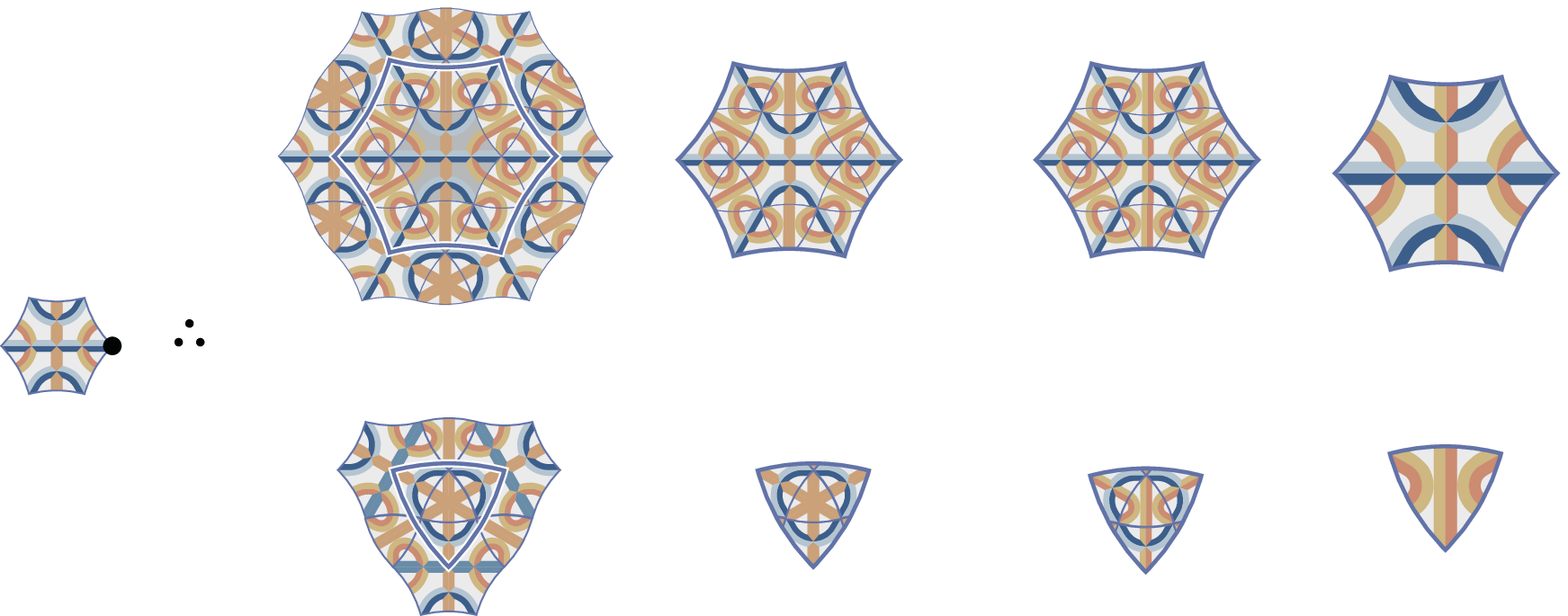}}

\vp Continuing this process, for illustration: 

{\includegraphics{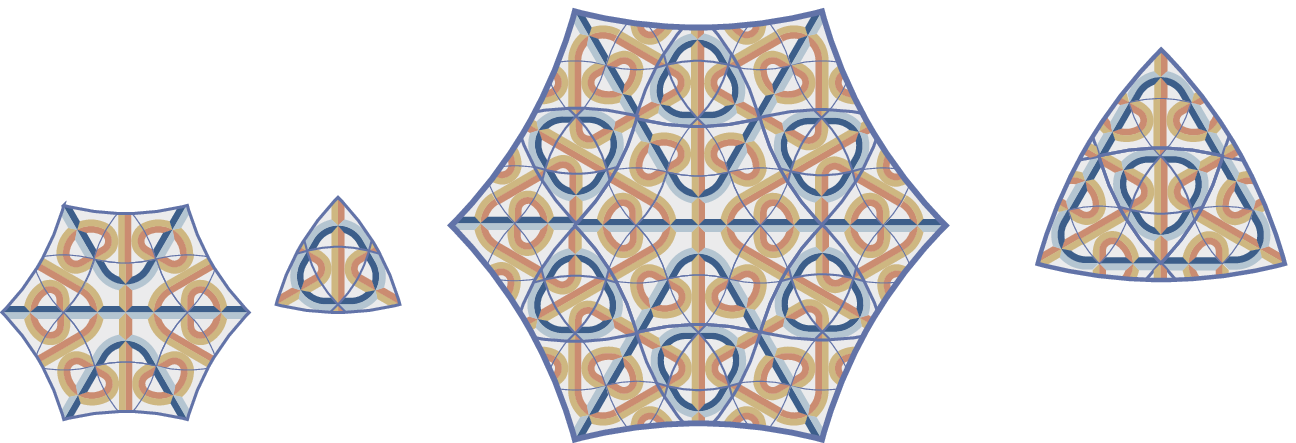}}

\vp A large patch:

\centerline{\includegraphics[width=1.1\textwidth]{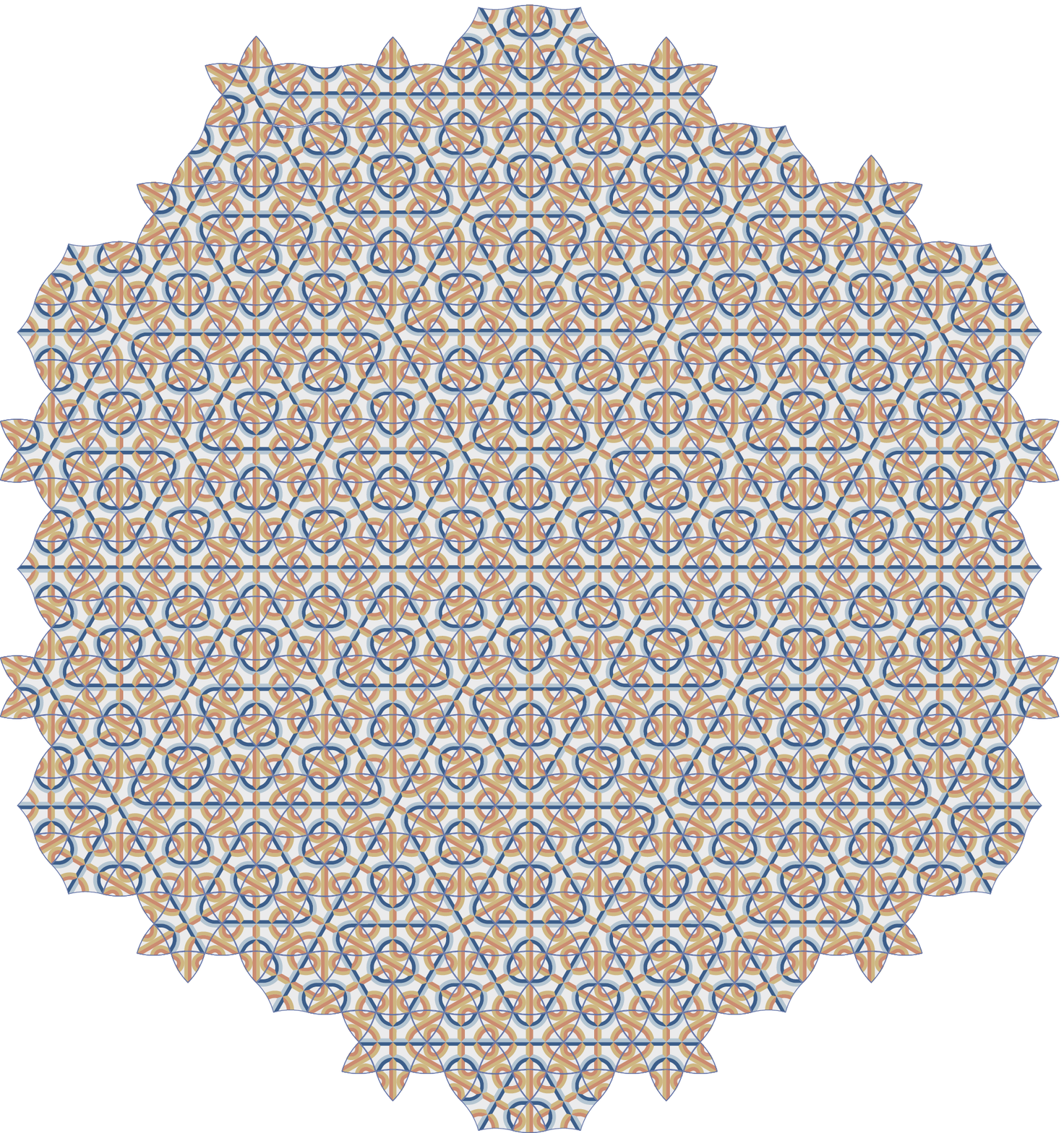}}

%




\end{document}